\documentclass[12pt]{article}%
\usepackage{amsmath}
\usepackage{amsfonts}
\usepackage{amssymb}
\usepackage{graphicx}
\providecommand{\U}[1]{\protect\rule{.1in}{.1in}}
\newtheorem{theorem}{Theorem}

\begin{document}

\title{Graham Higman's PORC theorem}
\author{Michael Vaughan-Lee}
\date{July 2018}
\maketitle

\begin{abstract}
Graham Higman published two important papers in 1960. In the first of these
papers he proved that for any positive integer $n$ the number of groups of
order $p^{n}$ is bounded by a polynomial in $p$, and he formulated his famous
PORC conjecture about the form of the function $f(p^{n})$ giving the number of
groups of order $p^{n}$. In the second of these two papers he proved that the
function giving the number of $p$-class two groups of order $p^{n}$ is PORC.
He established this result as a corollary to a very general result about
vector spaces acted on by the general linear group. This theorem takes over a
page to state, and is so general that it is hard to see what is going on.
Higman's proof of this general theorem contains several new ideas and is quite
hard to follow. However in the last few years several authors have developed
and implemented algorithms for computing Higman's PORC formulae in
special cases of his general theorem. These algorithms give perspective on
what are the key points in Higman's proof, and also simplify parts of the proof.

In this note I give a proof of Higman's general theorem written in the light
of these recent developments.

\end{abstract}

\section{Introduction}

Graham Higman wrote two immensely important and influential papers on
enumerating $p$-groups in the late 1950s. The papers were entitled
\emph{Enumerating }$p$\emph{-groups} I and II, and were published in the
Proceedings of the London Mathematical Society in 1960 (see \cite{higman60}
and \cite{higman60b}). In the first of these papers Higman proves that if we
let $f(p^{n})$ be the number of $p$-groups of order $p^{n}$, then%
\[
p^{\frac{2}{27}n^{2}(n-6)}\leq f(p^{n})\leq p^{(\frac{2}{15}+\varepsilon
_{n})n^{3}},
\]
where $\varepsilon_{n}$ tends to zero as $n$ tends to infinity. Higman also
formulated his famous PORC conjecture concerning the form of the function
$f(p^{n})$. He conjectured that for each $n$ there is an integer $N$
(depending on $n$) such that for $p$ in a fixed residue class modulo $N$ the
function $f(p^{n})$ is a polynomial in $p$. For example, {for }$p\geq5$ the
number of groups of order $p^{6}$ is{%
\[
3p^{2}+39p+344+24\gcd(p-1,3)+11\gcd(p-1,4)+2\gcd(p-1,5).
\]
(See \cite{newobvl}.) So for }$p\geq5$, $f(p^{6})$ is one of 8 polynomials in
$p$, with the choice of polynomial depending on the residue class of $p$
modulo 60. The number of groups of order $p^{6}$ is \textbf{P}olynomial
\textbf{O}n \textbf{R}esidue \textbf{C}lasses. In \cite{higman60b} Higman
proved that, for any given $n$, the function enumerating the number of
$p$-class 2 groups of order $p^{n}$ is a PORC function of $p$. He obtained
this result as a corollary to a very general theorem about vector spaces acted
on by the general linear group. As another corollary to this general theorem,
he also proved that for any given $n$ the function enumerating the number of
algebras of dimension $n$ over the field of $q$ elements is a PORC function of
$q$.

In recent years several authors have developed algorithms for computing
Higman's PORC formulae in various applications of his general theorem. Witty
\cite{witty} wrote a thesis describing an algorithm for computing the number
of $r$-generator $p$-class two groups, and I have published a series of papers
on this topic and on computing the numbers of non-associative algebras of
dimension $d$ (\cite{vlee12a}, \cite{vlee13}, \cite{vlee17}). \ Eick and
Wesche \cite{eickwesche} describe an algorithm for computing the numbers of
associative algebras over a finite field which are nilpotent of class $2$, and
this algorithm has been implemented in \textsf{GAP}. These algorithms simplify
parts of Higman's theory, and have given me a better understanding of his
general theorem and its proof. In this note I offer my insights into Higman's
remarkable theorem.

\section{Algebraic families of groups}

Higman introduces the notion of an \emph{algebraic family of groups}. Let
$\mathbb{Q}$ be the rational field and suppose we have a homomorphism%
\[
\varphi:\,\text{GL}(m,\mathbb{Q})\rightarrow\,\text{GL}(n,\mathbb{Q})
\]
with the property that if $A$ is a matrix in GL$(m,\mathbb{Q})$ then
$\varphi(A)$ is a matrix in GL$(n,\mathbb{Q})$ with entries of the form
$\frac{r}{s}$ where $r$ and $s$ are polynomials over $\mathbb{Q}$ in the
entries of $A$. (Of course the polynomials $r,s$ should depend only on
$\varphi$ and not on $A$.) Paraphrasing Higman slightly, he writes
\textquotedblleft Of course, the least common multiple of the denominators $s
$ is a power of $\det(A)$\textquotedblright. This stumped me for quite a
while, but eventually I consulted an expert in algebraic groups who assured me
that this was a well known, basic fact. In the end I managed to find my own
elementary proof of this, but I have decided not to include my proof in this
note. So we assume that the entries in 
$\varphi(A)$ are of the form $\frac{r}{\det(A)^k}$ where $r$ is a polynomial
over $\mathbb{Q}$ in the entries of $A$.

Higman's idea is that if $K$ is any
field whose characteristic does not divide the denominator of any of the
coefficients in the polynomials giving the entries in $\varphi(A)$ then the
coefficients in these polynomials can be interpreted as elements of the prime
subfield of $K$ so that $\varphi$ defines a homomorphism
\[
\varphi_{K}:\,\text{GL}(m,K)\rightarrow\,\text{GL}(n,K).
\]
Higman calls the collection of images $\varphi_{K}($GL$(m,K))$ an algebraic
family of groups.

Now let $\varphi:\,$GL$(m,\mathbb{Q})\rightarrow\,$GL$(n,\mathbb{Q})$ give an
algebraic family of groups, and let $P$ be the (finite) set of primes which
divide denominators of coefficients in the polynomials giving the
entries in $\varphi(A)$. Suppose that $K$ is a finite field of order $q$,
where the characteristic of $K$ is not contained in $P$. Let $V$ be a vector
space of dimension $n$ over $K$. Then $\varphi_{K}$ gives us an action of
GL$(m,K)$ on $V$. Higman \cite{higman60b} proves the following theorem.

\begin{theorem}
(a) The number of orbits of $V$ under the action of GL$(m,K)$, considered as a
function of $q$, is PORC.

(b) For each integer $k$ with $0\leq k\leq n$, the number of orbits of
GL$(m,K)$ on subspaces of $V$ of dimension $k$, considered as a function of
$q$, is PORC.
\end{theorem}

Actually Higman's algebraic families of groups are more general than this, as
they are given by homomorphisms%
\[
\varphi:\,\text{GL}(m_{1},\mathbb{Q})\times\,\text{GL}(m_{2},\mathbb{Q}%
)\times\ldots\times\,\text{GL}(m_{r},\mathbb{Q})\rightarrow\,\text{GL}%
(n,\mathbb{Q}).
\]
It may be that Higman's main reason for this generalization is that he proves
(b) (for $r=1$) by considering the action of GL$(k,K)\times\,$GL$(m,K)$ on the
space of $k\times n$ matrices, with GL$(k,K)$ acting on the left by matrix
multiplication, and GL$(m,K)$ acting on the right via $\varphi_{K}$. So to
prove (b) for $r=1$ he proves (a) for $r=2$, and more generally to prove (b)
for a direct product of $r$ general linear groups he proves (a) for a direct
product of $r+1$ general linear groups. The proof of Higman's theorem
given below avoids this complication, so I will restrict my attention to
the case $r=1$. Nevertheless the proof given here can
easily be adapted to prove Theorem 1 for Higman's more general algebraic
families, the only difference being that the notation would be more complicated.

\section{Two examples of algebraic families of groups}

\subsection{GL$(V)$ acting on $V\wedge V$}

Suppose that $V$ is a vector space of dimension $n$ over a field $K$. Then
there is a natural action of GL$(V)$ on the exterior square $V\wedge V$, which
we can describe as follows.

Suppose that $V$ has basis $v_{1},v_{2},\ldots,v_{n}$ over $K$, and let
$g\in\,$GL$(V)$ have matrix $A=\left[  a_{ij}\right]  $ with respect to this
basis, so that $v_{i}g=\sum_{j}a_{ij}v_{j}$. Then $V\wedge V$ has a basis
consisting of the elements $v_{i}\wedge v_{j}$ with $i<j$, and
\[
(v_{i}\wedge v_{j})g=\sum_{k,m}a_{ik}a_{jm}(v_{k}\wedge v_{m})=\sum
_{k<m}(a_{ik}a_{jm}-a_{im}a_{jk})(v_{k}\wedge v_{m}).
\]
So the matrix giving the action of $g$ on $V\wedge V$ with respect to the
basis $v_{i}\wedge v_{j}$ with $i<j$ has entries of the form $a_{ik}%
a_{jm}-a_{im}a_{jk}$. We have a homomorphism
\[
\varphi:\,GL(n,K)\rightarrow\,GL(\binom{n}{2},K)
\]
and if $A\in\,$GL$(n,K)$ then the entries in $\varphi(A)$ are integer
polynomials in the entries of $A$. So we have an algebraic family of groups.
Theorem 1 implies that if $K$ is a field of order $q$ then the number of
orbits of GL$(n,K)$ on $V\wedge V$ is PORC as a function of $q$, 
as is the number of orbits of GL$(n,K)$ on subspaces of $V\wedge V$ of 
dimension $k$ ($0\leq k\leq\binom{n}{2}$).

\subsection{Algebras over a field $K$}

Higman's general theorem also implies that for every dimension $m$ the number
of algebras of dimension $m$ over a field $K$ of order $q$ is a PORC function
of $q$. (See \cite{vlee13}.) By \textquotedblleft algebra\textquotedblright%
\ we mean a vector space with a bilinear product. There is no requirement that
the product satisfy any other condition such as associativity. If $B$ is an
algebra of dimension $m$ over $K$, and if we pick a basis $v_{1},v_{2}%
,\ldots,v_{m}$ for $B$ as a vector space over $K$ then for each pair of basis
elements $v_{i},v_{j}$ we can express the product $v_{i}v_{j}$ as a linear
combination%
\[
v_{i}v_{j}=\sum_{k}\lambda_{ijk}v_{k}%
\]
for some scalars $\lambda_{ijk}\in K$. These scalars are \emph{structure
constants} for the algebra $B$, and completely determine $B$. If we pick
another basis $w_{1},w_{2},\ldots,w_{m}$, and if%
\[
w_{i}w_{j}=\sum_{k}\mu_{ijk}w_{k}%
\]
then we obtain another set of structure constants $\mu_{ijk}$. We can express
the elements of the second basis as linear combinations of elements of the
first basis, and vice versa:%
\begin{align*}
w_{i}  & =\sum_{j=1}^{m}a_{ji}v_{j}\;(1\leq i\leq m),\\
v_{j}  & =\sum_{k=1}^{m}b_{kj}w_{k}\;(1\leq j\leq m),
\end{align*}
where $\left[  a_{ji}\right]  $ and $\left[  b_{kj}\right]  $ are $m\times m$
matrices over $K$ which are inverse to each other. So%
\begin{align*}
w_{i}w_{j}  & =\sum_{r,s=1}^{m}a_{ri}a_{sj}v_{r}v_{s}\\
& =\sum_{r,s,t=1}^{m}a_{ri}a_{sj}\lambda_{rst}v_{t}\\
& =\sum_{r,s,t,k=1}^{m}a_{ri}a_{sj}\lambda_{rst}b_{kt}w_{k}.
\end{align*}
It follows that%
\[
\mu_{ijk}=\sum_{r,s,t=1}^{n}a_{ri}a_{sj}\lambda_{rst}b_{kt}.
\]
If we think of the sets of structure constants as vectors in an $m^{3}$
dimensional vector space over $K$ then this gives us a homomorphism from
GL$(m,K)$ into GL$(m^{3},K)$ where the image of a matrix $A$ in GL$(m,K)$ has
entries of the form $\frac{f}{\det A}$ where $f$ is an integer polynomial in
the entries of $A$. Two sets of structure constants give isomorphic algebras
if and only if they lie in the same orbit under the action of GL\thinspace
$(m,K)$, and so Theorem 1 (a) implies that the number of $m$-dimensional
algebras over $K$, considered as a function of $q=|K|$, is PORC.

\section{Diagonal matrices in GL$(m,\mathbb{Q})$}

Let $\varphi:\,$GL$(m,\mathbb{Q})\rightarrow\,$GL$(n,\mathbb{Q})$ give an
algebraic family of groups. If $A\in \,$GL$(m,\mathbb{Q})$ then the entries
in $\varphi(A)$ have the form $\frac{r}{\det(A)^k}$ where $r$ is a polynomial
over $\mathbb{Q}$ in the entries of $A$, and we let $P$ be the finite set of primes 
which divide denominators of coefficients in the polynomials $r$. Let $R$ be 
the ring of rationals of the form $\frac{a}{b}$ where only primes in $P$ divide $b$.

\begin{theorem}
There is a matrix $Q$ with entries in $R$ and with $\det Q=\pm1$ such that if
$A$ is a diagonal matrix in GL$(m,\mathbb{Q})$ then $Q^{-1}\varphi(A)Q$ is
diagonal. Furthermore, if $A$ has eigenvalues $\lambda_{1},\lambda_{2}%
,\ldots,\lambda_{m}$ then $\varphi(A)$ has eigenvalues of the form
$\lambda_{1}^{n_{1}}\lambda_{2}^{n_{2}}\ldots\lambda_{m}^{n_{m}}$ for some
integers $n_{i}$.
\end{theorem}

To my mind this is the cleverest and trickiest part of Higman's proof of
Theorem 1. But note that in the case of the two examples given in Section 3
there is nothing to prove. In the first example $\varphi(A)$ is diagonal with
eigenvalues $\lambda_{i}\lambda_{j}$ ($i<j$), and in the second example
$\varphi(A)$ is diagonal with eigenvalues $\frac{\lambda_{i}\lambda_{j}%
}{\lambda_{k}}$ ($i,j,k=1,2,\ldots,m$). Similarly in Eick and Wesche's
algorithm \cite{eickwesche} to compute the numbers of class two associative
algebras they consider the action of GL$(V)$ on $V\otimes V$, and in this case
if $A$ is diagonal with eigenvalues $\lambda_{1},\lambda_{2},\ldots
,\lambda_{m}$ then $\varphi(A)$ is diagonal with eigenvalues $\lambda
_{i}\lambda_{j}$ ($i,j=1,2,\ldots,m$). Of course this is only true if we
choose the \textquotedblleft right\textquotedblright\ basis for $V\otimes V$.
Theorem 2 implies that if we have an algebraic family of groups giving an action 
of GL$(m,\mathbb{Q})$ on a vector space $W$, then we can always choose a basis 
of $W$ with respect to which diagonal matrices in GL$(m,\mathbb{Q})$ act 
diagonally on $W$.

Let $A\in\,$GL$(m,\mathbb{Q})$ be a diagonal matrix with eigenvalues
$\lambda_{1},\lambda_{2},\ldots,\lambda_{m}$. The entries in $\varphi(A)$ 
are $R$-linear combinations of products
$\lambda_{1}^{n_{1}}\lambda_{2}^{n_{2}}\ldots\lambda_{m}^{n_{m}} $ with
$n_{i}\in\mathbb{Z}$. Let $\lambda_{1}^{n_{i1}%
}\lambda_{2}^{n_{i2}}\ldots\lambda_{m}^{n_{im}}$ ($1\leq i\leq k$) be the
distinct products of eigenvalues of $A$ and their inverses which occur in
$\varphi(A)$. We can write%
\[
\varphi(A)=\sum_{i=1}^{k}E_{i}\lambda_{1}^{n_{i1}}\lambda_{2}^{n_{i2}}%
\ldots\lambda_{m}^{n_{im}},
\]
where $E_{1},E_{2},\ldots,E_{k}$ are $n\times n$ matrices with entries in $R
$. If we let $B$ be a diagonal matrix in GL$(m,\mathbb{Q)}$ with eigenvalues
$\mu_{1},\mu_{2},\ldots,\mu_{m}$ then%
\[
\varphi(B)=\sum_{i=1}^{k}E_{i}\mu_{1}^{n_{i1}}\mu_{2}^{n_{i2}}\ldots\mu
_{m}^{n_{im}}%
\]
and%
\[
\varphi(AB)=\sum_{i=1}^{k}E_{i}(\lambda_{1}\mu_{1})^{n_{i1}}(\lambda_{2}%
\mu_{2})^{n_{i2}}\ldots(\lambda_{m}\mu_{m})^{n_{im}}.
\]
Since $\varphi(A)\varphi(B)=\varphi(AB)$ for all $\lambda_{1},\lambda
_{2},\ldots,\lambda_{m},\mu_{1},\mu_{2},\ldots,\mu_{m}\in\mathbb{Q}\setminus \{0\}$ this
implies that $E_{i}E_{j}=0$ for all $i\neq j$ and that $E_{i}^{2}=E_{i}$ for
all $i$.

So the matrices $E_{1},E_{2},\ldots,E_{k}$ can be simultaneously diagonalized.
This means that we can find a non-singular $n\times n$ matrix $C$ such that
every column of $C$ is an eigenvector with eigenvalue $0$ or $1$ for each of
$E_{1},E_{2},\ldots,E_{k}$. We can take the entries in $C$ to be integers.
Since $E_{i}E_{j}=0$ for $i\neq j$ and since $\varphi(A)$ is non-singular,
each column of $C$ has eigenvalue 1 for exactly one of the matrices
$E_{1},E_{2},\ldots,E_{k}$. We can order the columns of $C$ so that the first
few columns are eigenvectors with eigenvalue 1 for $E_{1}$, so that the next
few columns are eigenvectors with eigenvalue 1 for $E_{2}$, and so on. For
$i=1,2,\ldots,n$ let $e_{i}$ be the column vector with $1$ in the $i^{th}$
place, and $0$ in every other place.

We consider elementary row operations on $C$ of the following three forms:

\begin{enumerate}
\item Swap two rows.

\item Subtract an integer multiple of one row from another.

\item Multiply a row by $-1$.
\end{enumerate}

If $a$ is the greatest common divisor of the entries in the first column of
$C$ then we can apply a sequence of row operations to reduce the first column
of $C$ to $ae_{1}$. Then we can apply a sequence of row operations to rows
$2,3,\ldots,n$ to reduce the second column of $C$ to $be_{1}+ce_{2}$ for some
$b,c$. Next we apply a sequence of row operations to rows $3,4,\ldots,n $ to
reduce the third column to $de_{1}+ee_{2}+fe_{3}$ for some $d,e,f$. Continuing
in this way we eventually reduce $C$ to an upper triangular integer matrix.
Applying this sequence of row operations to $C$ corresponds to premultiplying
$C$ by a sequence of elementary matrices. Multiplying these elementary
matrices together we obtain an integer matrix $Q$ with $\det Q=\pm1$ such that
$QC$ is upper triangular. Let $QC=D$ and let $F_{i}=QE_{i}Q^{-1}$ for
$i=1,2,\ldots,k$. Then $F_{i}F_{j}=0$ for $i\neq j$ and $F_{i}^{2}=F_{i}$ for
all $i$. Let $V_{i}=\ker(F_{i}-1)$ and let $\dim V_{i}=d_{i}$ for
$i=1,2,\ldots,k$. Then the first $d_{1}$ columns of $D$ form a basis for
$V_{1}$, the next $d_{2}$ columns of $D$ form a basis for $V_{2}$, and so on. However

\begin{itemize}
\item the first $d_{1}$ columns of $D$ span the same space as $e_{1}%
,e_{2},\ldots,e_{d_{1}}$,

\item the next $d_{2}$ columns span the same space as%
\[
e_{d_{1}+1}+v_{1},e_{d_{1}+2}+v_{2},\ldots,e_{d_{1}+d_{2}}+v_{d_{2}}%
\]
for some $v_{1},v_{2},\ldots,v_{d_{2}}\in V_{1},$

\item the next $d_{3}$ columns span the same space as%
\[
e_{d_{1}+d_{2}+1}+w_{1},e_{d_{1}+d_{2}+2}+w_{2},\ldots,e_{d_{1}+d_{2}+d_{3}%
}+w_{d_{3}}%
\]
for some $w_{1},w_{2},\ldots,w_{d_{3}}\in V_{1}+V_{2},$

\item and so on.
\end{itemize}

Now let $E$ be the $n\times n$ matrix with columns%
\[
e_{1},e_{2},\ldots,e_{d_{1}},e_{d_{1}+1}+v_{1},\ldots,e_{d_{1}+d_{2}}%
+v_{d_{2}},e_{d_{1}+d_{2}+1}+w_{1},\ldots,e_{d_{1}+d_{2}+d_{3}}+w_{d_{3}%
},\ldots,
\]
so that the first $d_{1}$ columns of $E$ form a basis for $V_{1}$, the next
$d_{2}$ columns form a basis for $V_{2}$, and so on. Then%
\[
e_{d_{1}+1}+v_{1}=F_{2}(e_{d_{1}+1}+v_{1})=F_{2}e_{d_{1}+1}%
\]
since $F_{1}F_{2}=0$. All the entries in the matrix $F_{2}$ lie in $R$, and so
all the entries in $v_{1}$ lie in $R$. Similarly, all the entries in
$v_{2},v_{3},\ldots,v_{d_{2}},w_{1},\ldots,w_{d_{3}},\ldots$ lie in $R$. So
$E$ is an upper triangular matrix with 1's down the diagonal and with all the
entries above the diagonal lying in the ring $R$. It follows that $Q^{-1}E$ is
a matrix with entries in $R$ and with determinant $\pm1$, and such that
$E^{-1}Q\varphi(A)Q^{-1}E$ is diagonal. This completes the proof of Theorem 2.

From now on we replace $\varphi:\,$GL$(m,\mathbb{Q})\rightarrow\,$%
GL$(n,\mathbb{Q})$ by $\varphi^{\ast}$, where
\[
\varphi^{\ast}(A)=E^{-1}Q\varphi(A)Q^{-1}E.
\]
In other words, we assume that $\varphi(A)$ is diagonal whenever $A$ is diagonal.

\section{Matrices in Jordan form}

For each integer $k\geq1$ let $J_{k}$ be the $k\times k$ Jordan matrix with
1's down the diagonal and 1's down the superdiagonal. Higman considers a
non-singular matrix to be in Jordan form if it can be expressed in the form%
\begin{equation}
\lambda_{1}J_{k_{1}}\oplus\lambda_{2}J_{k_{2}}\oplus\ldots\oplus\lambda
_{r}J_{k_{r}}%
\end{equation}
for some eigenvalues $\lambda_{1},\lambda_{2},\ldots,\lambda_{r}$ and some
integers $k_{1},k_{2},\ldots,k_{r}$. (This is possible since the eigenvalues
are non-zero.) So let $A\in\,$GL$(m,\mathbb{Q})$ be a matrix of form (1). For
the moment assume that the eigenvalues $\lambda_{1},\lambda_{2},\ldots
,\lambda_{r}$ are all distinct. Actually, it helps to think of $\lambda
_{1},\lambda_{2},\ldots,\lambda_{r}$ as indeterminates. Let%
\begin{gather*}
\Lambda=\lambda_{1}I_{k_{1}}\oplus\lambda_{2}I_{k_{2}}\oplus\ldots
\oplus\lambda_{r}I_{k_{r}},\\
J=J_{k_{1}}\oplus J_{k_{2}}\oplus\ldots\oplus J_{k_{r}}.
\end{gather*}

Then $A=\Lambda J=J\Lambda$. By Theorem 2 we may suppose that $\varphi
(\Lambda)$ is diagonal, with eigenvalues which are products of the eigenvalues
$\lambda_{1},\lambda_{2},\ldots,\lambda_{r}$ and their inverses. Suppose that
the products that arise as eigenvalues of $\varphi(\Lambda)$ are $m_{1}%
,m_{2},\ldots,m_{s}$, and let $Q$ be a permutation matrix chosen so that%
\[
L=Q^{-1}\varphi(\Lambda)Q=m_{1}I_{t_{1}}\oplus m_{2}I_{t_{2}}\oplus
\ldots\oplus m_{s}I_{t_{s}}%
\]
for some positive integers $t_{1},t_{2},\ldots,t_{s}$. Let $M=Q^{-1}%
\varphi(J)Q$. Then $LM=ML$ and so%
\[
M=E_{1}\oplus E_{2}\oplus\ldots\oplus E_{s}%
\]
where $E_{i}$ is a $t_{i}\times t_{i}$ matrix for $i=1,2,\ldots,s$, and
where the entries in $E_i$ all lie in $R$ . The
matrix $J$ is conjugate to all its power $J^{i}$ ($i=1,2,\ldots$) and so $M$
is also conjugate to all its powers. This implies that 1 is the only
eigenvalue of $M$. And this implies that we can find invertible matrices
$X_{1},X_{2},\ldots,X_{s}$ with rational entries such that $X_{i}^{-1}E_{i}X_{i}$ is 
in Jordan form (with 1 as the only eigenvalue) for $i=1,2,\ldots,s$. Let%
\[
X=X_{1}\oplus X_{2}\oplus\ldots\oplus X_{s}.
\]
So%
\begin{equation}
X^{-1}Q^{-1}\varphi(A)QX=LX^{-1}MX
\end{equation}
is in Jordan form. Note that the matrices $Q,X$ do not depend on the values of
$\lambda_{1},\lambda_{2},\ldots,\lambda_{r}$. Furthermore $LX^{-1}MX$ is in
Jordan form even if $\lambda_{1},\lambda_{2},\ldots,\lambda_{r}$ are not all
distinct. Now suppose that $p$ is a prime which does not divide the
denominator of any of the entries in $X$ and does not divide the numerator of
$\det X$. If $K$ is a finite field of characteristic $p$ then we can interpret
$X$ and $Q$ as non singular matrices with entries in the prime subfield of
$K$, so equation (2) also gives the Jordan form of $\varphi_{K}(A)$ if $A$ is
a matrix in GL$(m,K)$ of the form (1). Note that the sizes of the Jordan
blocks in $X^{-1}MX$ and the number of blocks of each size depends only on the
integers $k_{1},k_{2},\ldots,k_{r}$. Also the eigenvalues of $\varphi_{K}(A)$
corresponding to the Jordan blocks in $X^{-1}MX$ have the form
\[
m_i=\lambda_{1}^{n_{i1}}\lambda_{2}^{n_{i2}}\ldots \lambda_{r}^{n_{ir}}
\]
for some integers $n_{ij}$ which depend only on $k_{1},k_{2},\ldots,k_{r}$.

There will be a finite
number of \textquotedblleft exceptional\textquotedblright\ characteristics
which divide one of the denominators of the entries in $X$, or divide the
numerator of $\det X$. Let $K$ be a finite field with exceptional
characteristic $p$. Let $A$ be a matrix in GL$(m,K)$ of the form (1). We
follow the same analysis as above and obtain the same expression%
\[
Q^{-1}\varphi_{K}Q=LM=ML
\]
as above. Now $J^{p^{m}}=I_{m}$, and this implies that $\varphi_{K}(J^{p^{m}%
})=I_{n}$. If $\lambda$ is an eigenvalue of $\varphi_{K}(J)$ then
$\lambda^{p^{m}}$ is an eigenvalue of $\varphi_{K}(J^{p^{m}})$, and so
$\lambda^{p^{m}}=1$ which implies that $\lambda=1$. So, just as above, we can
find invertible matrices $X_{1},X_{2},\ldots,X_{s}$ with entries in GF$(p)$ such 
that $X_{i}^{-1}E_{i}X_{i}$ is in Jordan form (with 1 as the only eigenvalue) for
$i=1,2,\ldots,s$. This gives the Jordan form of $\varphi_{K}(A)$ for all
fields of characteristic $p$.

\bigskip
It might be helpful to give a simple example. Let $K$ be a field and let
$A=aJ_{2}\oplus bJ_{3}$ for some $a,b\in K$. If the characteristic of $K$ is 0
or is a prime $p>3$ then the Jordan form of $A\otimes A$ is%
\[
a^{2}J_{1}\oplus a^{2}J_{3}\oplus abJ_{2}\oplus abJ_{2}\oplus abJ_{4}\oplus
abJ_{4}\oplus b^{2}J_{1}\oplus b^{2}J_{3}\oplus b^{2}J_{5}.
\]
The exceptional characteristics are $p=2,3$. In characteristic 2 the Jordan
form of $A\otimes A$ is%
\[
a^{2}J_{2}\oplus a^{2}J_{2}\oplus abJ_{2}\oplus abJ_{2}\oplus abJ_{4}\oplus
abJ_{4}\oplus b^{2}J_{1}\oplus b^{2}J_{4}\oplus b^{2}J_{4},
\]
and in characteristic 3 it is%
\[
a^{2}J_{1}\oplus a^{2}J_{3}\oplus abJ_{3}\oplus abJ_{3}\oplus abJ_{3}\oplus
abJ_{3}\oplus b^{2}J_{3}\oplus b^{2}J_{3}\oplus b^{2}J_{3}.
\]
This example illustrates that it is much easier in practice than in theory to
show for a given characteristic that when $A$ is of form (1) then the Jordan
form of $\varphi(A)$ depends only on the integers $k_{1},k_{2},\ldots,k_{r}$
and on the values of $\lambda_{1},\lambda_{2},\ldots,\lambda_{r}$. There are
only finitely many possible choices for $k_{1},k_{2},\ldots,k_{r}$ which sum
to $m$. For each such choice we treat the eigenvalues $\lambda_{1},\lambda
_{2},\ldots,\lambda_{r}$ as indeterminates and compute the Jordan form of
$\varphi(A)$ in characteristic zero. We then identify the exceptional
characteristics, and compute the Jordan form in each exceptional characteristic.

Note that this example also covers the case when $a=b$ (in which case
$A\otimes A$ has only 1 eigenvalue) and the case $a=-b$ (when $A\otimes A$ has
two eigenvalues). In all other cases $A\otimes A$ has three eigenvalues.

\section{The type of a matrix in GL$(m,K)$}

Let $K$ be a field, and let $A$ be a matrix in GL$(m,K)$. Let the primary
invariant factors of $A$ be $p_{1}(x)^{e_{1}},p_{2}(x)^{e_{2}},\ldots
,p_{k}(x)^{e_{k}}$ where $p_{1},p_{2},\ldots,p_{k}$ are monic irreducible
polynomials in $K[x]$. Let the distinct irreducible polynomials which occur in
the sequence $p_{1},p_{2},\ldots,p_{k}$ be $q_{1},q_{2},\ldots,q_{t}$ (with
$t\leq k$). For $1\leq i\leq t$ let $S_{i}$ be the multiset of exponents $e$
such that $q_{i}^{e}$ is one of the primary invariant factors of $A$. Then the
type of $A$ is the multiset of ordered pairs%
\[
\{(\deg q_{1},S_{1}),(\deg q_{2},S_{2}),\ldots,(\deg q_{t},S_{t})\}.
\]
For example, if the primary invariant factors of $A$ are $p(x)^{2}$,
$p(x)^{3}$, $q(x)$, $q(x)$, $q(x)^{4}$ where $p$ and $q$ are distinct monic 
irreducible polynomials, then the type of $A$ is%
\[
\{(\deg p,\{2,3\}),(\deg q,\{1,1,4\})\}.
\]
(Note that repeated entries in these multisets are significant.) So the type
of $A$ records the degrees of the different irreducible polynomials which
arise in the primary invariant factors of $A$, together with the multiset of
exponents associated with each of these irreducible polynomials. There are
only finitely many possible types of matrices in GL$(m,K)$. In addition if $K
$ has order $q$ then the number of matrices in GL$(m,K)$ of a given type is a
polynomial in $q$. Green \cite{green55} proves that the size of the conjugacy
class of $A$ is a polynomial in $q$, with the polynomial depending only on the
type of $A$. A formula for this polynomial is given on page 181 of
\cite{macdonald79}.

If $A$ is a matrix in Jordan form, and if $A$ has the form (1) from
Section 5, then the type of $A$ depends only on which equalities
$\lambda_{i}=\lambda_{j}$ ($i<j$) hold between the eigenvalues of $A$. And for
a given characteristic the type of $\varphi(A)$ depends only on which equalities
$m_{i}=m_{j}$ ($i<j$) hold between the eigenvalues of $\varphi(A) $. Thus in
the example given at the end of the last section the type of $A$ depends only
on whether $a=b$ or not, and the type of $A\otimes A$ depends on whether $a=b$
or $a=-b$.

We want to use the results of Section 5 to compute the type of $\varphi(A)$
even when the characteristic polynomial of $A$ does not split into linear
factors, and we proceed as follows. Suppose that $A\in\,$GL$(m,K)$ where $K$
is a finite field of order $q$. Let $L$ be the splitting field of the
characteristic polynomial of $A$. Then let $B$ be the Jordan form of $A$ when
considered as a matrix in GL$(m,L)$. Suppose that $p(x)^{e}$ is a primary
invariant factor of $A$, where $p(x)$ has degree $d$, and where $\lambda$ is a
root of $p(x)$ in $L$. Then the rational canonical form of $A$ has the
companion matrix of $p(x)^{e}$ as one of its blocks, and corresponding to this
we have%
\[
\lambda J_{e}\oplus\lambda^{q}J_{e}\oplus\lambda^{q^{2}}J_{e}\oplus
\ldots\oplus\lambda^{q^{d-1}}J_{e}%
\]
as a sum of blocks in $B$. Using the results of Section 5, the Jordan form of
$\varphi_{L}(B)$ has the form%
\[
m_{1}J_{t_{1}}\oplus m_{2}J_{t_{2}}\oplus\ldots\oplus m_{s}J_{t_{s}}%
\]
for some positive integers $t_{1},t_{2},\ldots,t_{s}$, and some products
$m_{1},m_{2},\ldots,m_{s}$ of the eigenvalues of $B$ and their inverses. The
integers $t_{1},t_{2},\ldots,t_{s}$ depend only on the type of $A$ and the
characteristic of $K$.

We need to investigate the eigenvalues of $B$ and the products $m_{1}%
,\ldots,m_{s}$ more closely. Let the distinct irreducible
polynomials which divide the primary invariant factors of $A$ be $q_{1}%
,q_{2},\ldots,q_{t}$, and let $A$ have type%
\[
T=\{(n_{1},S_{1}),(n_{2},S_{2}),\ldots,(n_{t},S_{t})\},
\]
where $n_{i}=\deg q_{i}$. We pick a root $\lambda_{i}$ of the polynomial
$q_{i}$ for $i=1,2,\ldots,t$. Then the eigenvalues of $B$ are $\{\lambda
_{i}^{q^{r}}:\,1\leq i\leq t,\;0\leq r<n_{i}\}$, and%
\[
B=%
{\displaystyle\bigoplus\limits_{i=1}^{t}}
{\displaystyle\bigoplus\limits_{j\in S_{i}}}
{\displaystyle\bigoplus\limits_{k=0}^{n_{i}-1}}
\lambda_{i}^{q^k}J_{j}.
\]

The eigenvalues $\lambda_{i}$ satisfy the equations%
\begin{equation}
\lambda_{i}^{q^{n_{i}}}=\lambda_{i}\;(i=1,2,\ldots,t).%
\end{equation}
They also satisfy the non-equations%
\begin{equation}
\lambda_{i}^{q^{r}}\neq\lambda_{i}\;(0<r<n_{i},\;i=1,2,\ldots,t),%
\end{equation}
\begin{equation}
\lambda_{i}^{q^{r}}\neq\lambda_{j}^{q^{s}}\;(i\neq j,\;0\leq r<n_{i},\;0\leq s<n_{j}).%
\end{equation}
The products $m_{1},m_{2},\ldots,m_{s}$ giving the eigenvalues of $\varphi_{L}(B)$
are of the form%
\[
\lambda_{1}^{h_{1}}\lambda_{2}^{h_{2}}\ldots\lambda_{t}^{h_{t}}%
\]
where $h_{1},h_{2},\ldots,h_{t}$ are integer polynomials in $q$. Note that
these polynomials depend only on the type of $A$ and on the characteristic of
$K$. The map $m_{i}\longmapsto m_{i}^{q}$ gives a permutation of $m_{1}%
,m_{2},\ldots,m_{s}$ and we can work out this permutation using the relations
(3). Provided we know which relations $m_{i}=m_{j}$ hold for these particular
values of $\lambda_{1},\lambda_{2},\ldots,\lambda_{t}$, and also know which of
these relations do not hold, then we can work out the type of $\varphi_{K}%
(A)$. As we will see in Section 8, as well as being able to calculate the type
of $\varphi_{K}(A)$ we also need to be able to calculate the dimension of the
eigenspace of $\varphi_{K}(A)$ with eigenvalue 1. This is just the number of
eigenvalues $m_{i}$ ($1\leq i\leq s$) which are equal to 1 (for these particular
values of $\lambda_{1},\lambda_{2},\ldots,\lambda_{t}$).

Now suppose that $\mu_{1},\mu_{2},\ldots,\mu_{t}\in L$ satisfy the relations
(3) and the non-relations (4), (5). If we let $g_{i}$ be the minimum polynomial of
$\mu_{i}$ over $K$ for $i=1,2,\ldots,t$ then%
\[%
{\displaystyle\bigoplus\limits_{i=1}^{t}}
{\displaystyle\bigoplus\limits_{j\in S_{i}}}
{\displaystyle\bigoplus\limits_{k=0}^{n_{i}-1}}
\mu_{i}^{q^k}J_{j}%
\]
is conjugate to matrices in GL$(m,K)$ with primary invariant factors%
\[
\{g_{i}^{j}\,:\,1\leq i\leq t,\;j\in S_{i}\},
\]
and all these matrices have type $T$. In this sense, $\mu_{1},\mu_{2}%
,\ldots,\mu_{t}$ determines a conjugacy class of matrices of type $T$ in
GL$(m,K)$. As we range over all possible solutions $\mu_{1},\mu_{2},\ldots
,\mu_{t}$ the conjugacy classes determined by $\mu_{1},\mu_{2},\ldots,\mu_{t}$
range over all possible conjugacy classes of matrices of type $T$. Furthermore
each such conjugacy class arises the same number of times. We get the same
conjugacy class if we replace $\mu_{1},\mu_{2},\ldots,\mu_{t}$ by $\nu_{1}%
,\nu_{2},\ldots,\nu_{t}$ where $\nu_{i}$ is conjugate to $\mu_{i}$ for all
$i$. Also if $(n_{i},S_{i})=(n_{j},S_{j})$ then we obtain the same conjugacy
class if we swap $\mu_{i}$ and $\mu_{j}$. We can make this precise as follows.
Write the entries $(n_{i},S_{i})$ from $T$ in a list%
\[
\lbrack(n_{1},S_{1}),(n_{2},S_{2}),\ldots,(n_{t},S_{t})]
\]
and let $G$ be the group of permutations $\pi$ of $\{1,2,\ldots,t\}$ such that%
\[
\lbrack(n_{1},S_{1}),(n_{2},S_{2}),\ldots,(n_{t},S_{t})]=[(n_{1\pi},S_{1\pi
}),(n_{2\pi},S_{2\pi}),\ldots,(n_{t\pi},S_{t\pi})].
\]
Then as we range over all possible solutions in $L$ of (3), (4) and (5) we run
through all possible conjugacy classes of elements in GL$(m,K)$ with type $T$,
and each conjugacy class arises
\begin{equation}
n_{1}n_{2}\ldots n_{t}|G|
\end{equation}
times.
\bigskip

To help clarify these ideas we investigate two simple examples. Let $K$ be a
finite field of order $q$, and let $A\in\,$GL$(m,K)$ have primary
invariant factors $g^{2}h^{3}$ where $g$ and $h$ are different monic
irreducible polynomials of degree 2. So $A$ has type%
\[
\{(2,\{2\}),(2,\{3\})\}.
\]
Let $\lambda$ be a root of $g$ and let $\mu$ be a root of $h$ in the splitting
field $L$ of $gh$ over $K$. Then the Jordan form of $A$ over $L$ is%
\[
B=\lambda J_{2}\oplus\lambda^{q}J_{2}\oplus\mu J_{3}\oplus\mu^{q}J_{3}.
\]
The eigenvalues of $B$ satisfy $\lambda^{q^{2}}=\lambda$, $\mu^{q^{2}}=\mu$,
$\lambda\neq\lambda^{q}$, $\mu\neq\mu^{q}$, $\lambda\neq\mu$, $\lambda\neq
\mu^{q}$. These equalities and inequalities determine the type of $A$. (We
also have inequalities $\lambda^{q}\neq\mu$, $\lambda^{q}\neq\mu^{q}$ but these
are redundant.) We can run over all conjugacy classes of GL$(m,K)$ of type
$\{(2,\{2\}),(2,\{3\})\}$ by running over all possible choices of $\lambda
,\mu$ in GF$(q^{2})$ satisfying these equations and non equations, and each
conjugacy class will arise 4 times, since swapping $\lambda$ and $\lambda^{q}$
or $\mu$ and $\mu^{q}$ gives the same conjugacy class.

As a second example, suppose that $K$ is a finite field of order $q$, and that
$A\in\,$GL$(m,K)$ has a single primary invariant factor $g^{3}$ where $g$ is
an irreducible quadratic. So $A$ has type $\{2,\{3\}\}$. If we pick a root
$\lambda$ of $g$ in GF$(q^{2})$ then the Jordan form of $A$ when considered as
a matrix over GF$(q^{2})$ is $\lambda J_{3}\oplus\lambda^{q}J_{3}$. The
eigenvalue $\lambda$ satisfies $\lambda^{q^{2}}=\lambda$, $\lambda\neq
\lambda^{q}$. Provided the characteristic of $K$ is at least 5, the Jordan
form of $A\otimes A$ over GF$(q^{2})$ is%
\[
\lambda^{2}J_{1}\oplus\lambda^{2}J_{3}\oplus\lambda^{2}J_{5}\oplus
\lambda^{q+1}J_{1}\oplus\lambda^{q+1}J_{1}\oplus\lambda^{q+1}J_{3}%
\]
\[
\oplus\lambda^{q+1}J_{3}\oplus\lambda^{q+1}J_{5}\oplus\lambda^{q+1}J_{5}%
\oplus\lambda^{2q}J_{1}\oplus\lambda^{2q}J_{3}\oplus\lambda^{2q}J_{5}.
\]
All the eigenvalues of $A\otimes A$ lie in the set $\{\lambda^{2}%
,\lambda^{q+1},\lambda^{2q}\}$, and we have%
\[
(\lambda^{2})^{q}=\lambda^{2q},\;(\lambda^{q+1})^{q}=\lambda^{q+1}%
,\;(\lambda^{2q})^{q}=\lambda^{2}.
\]
We also have%
\[
\lambda^{2}\neq\lambda^{q+1},\;\lambda^{2q}\neq\lambda^{q+1},
\]
so the type of $A\otimes A$ as a matrix over $K$ depends on whether or not the
equation $\lambda^{2}=\lambda^{2q}$ is satisfied. So to compute the numbers of
times matrices $A\otimes A$ of these two types arise as $A$ ranges over
conjugacy classes of type $\{2,\{3\}\}$ we need to count the numbers of
choices of $\lambda$ in GF$(q^{2})$ which satisfy the following two sets of
equations and non-equations%
\begin{align*}
\lambda^{q^{2}}  & =\lambda,\;\lambda\neq\lambda^{q},\;\lambda^{2}%
=\lambda^{2q},\\
\lambda^{q^{2}}  & =\lambda,\;\lambda\neq\lambda^{q},\;\lambda^{2}\neq
\lambda^{2q}.
\end{align*}
(We need to divide these answers by 2 to account for the fact that swapping
$\lambda$ and $\lambda^{q}$ gives the same conjugacy class.) The dimension of
the eigenspace of $A\otimes A$ with eigenvalue 1 is six if $\lambda^{q+1}=1$,
and zero otherwise.

\section{Choosing elements from finite fields}

Higman \cite{higman60b} proves the following theorem.

\begin{theorem}
The number of ways of choosing a finite number of elements from GF$(q^{n})$
subject to a finite number of monomial equations and inequalities between them
and their conjugates over GF$(q)$, considered as a function of $q$, is PORC.
\end{theorem}

Here we are choosing elements $x_{1},x_{2},\ldots,x_{k}$ (say) from the finite
field GF$(q^{n})$ (where $q$ is a prime power) subject to a finite set of
equations and non-equations of the form%
\[
x_{1}^{n_{1}}x_{2}^{n_{2}}\ldots x_{k}^{n_{k}}=1
\]
and%
\[
x_{1}^{n_{1}}x_{2}^{n_{2}}\ldots x_{k}^{n_{k}}\neq1
\]
where $n_{1},n_{2},\ldots,n_{k}$ are integer polynomials in the Frobenius
automorphism $x\rightarrow x^{q}$ of GF$(q^{n})$. Higman calls these equations
and non-equations monomial. For example, as I showed in the first example at
the end of Section 6, one way of computing the number of conjugacy classes of
matrices $A\in\,$GF$(q)$ of type $\{(2,\{2\}),(2,\{3\})\}$ is to count the
number of choices of $\lambda,\mu$ in GF$(q^{2})$ satisfying%
\[
\lambda^{q^{2}}=\lambda,\,\mu^{q^{2}}=\mu,\,\lambda\neq\lambda^{q},\,\mu
\neq\mu^{q},\,\lambda\neq\mu,\,\lambda\neq\mu^{q},
\]
and then divide by 4. Of course you can write these equations and
non-equations as%
\[
\lambda^{q^{2}-1}=1,\,\mu^{q^{2}-1}=1,\,\lambda^{q-1}\neq1,\,\mu^{q-1}%
\neq1,\,\lambda\mu^{-1}\neq1,\,\lambda\mu^{-q}\neq1,
\]
to match Higman's notation. Higman's proof of Theorem 3 involves 5 pages of
homological algebra, but a shorter more elementary proof can be found in
\cite{vlee12a} and in \cite{vlee17}.

To prove Theorem 3 you actually only need to prove that the number of ways of
choosing a finite number of elements from GF$(q^{n})$ subject to a finite
number of monomial equations between them and their conjugates over GF$(q)$,
considered as a function of $q$, is PORC. To see this suppose that we have a
set $S$ of equations and a set $T$ of non-equations. Let $T^{\ast}$ be the set
of equations obtained from $T$ be replacing all the $\neq$'s by $=$'s. For
each subset $U\subseteq T^{\ast}$ let $n_{U}$ be the number of solutions to
the equations $S\cup U$. Then the number of solutions to the equations $S$ and
the non-equations $T$ is%
\[
\sum_{U\subseteq T^{\ast}}(-1)^{|U|}n_{U}\text{.}%
\]

In \cite{vlee12a} and in \cite{vlee17} I show that to find the number of ways
of choosing a finite number of elements from GF$(q^{n})$ subject to a finite
number of monomial equations $S$ we write the equations in $S$ as the rows of
a matrix. We also have to add in equations $x_{i}^{q^{n}-1}=1$ to make sure
that the solutions lie in GF$(q^{n})$. For example, we represent the equations%
\[
x_{1}^{q^{2}-1}=1,\;x_{1}^{q+1}x_{2}^{-2}=1,\;x_{1}^{q^{n}-1}=1,\;x_{2}%
^{q^{n}-1}=1
\]
by the matrix%
\[
\left[
\begin{array}
[c]{cc}%
q^{2}-1 & 0\\
q+1 & -2\\
q^{n}-1 & 0\\
0 & q^{n}-1
\end{array}
\right]  .
\]
For any given value of $q$ this matrix is an integer matrix and the number of
solutions to the equations is the product of the elementary divisors in the
Smith normal form of the matrix. In \cite{vlee17} I show that the the number
of solutions to a set of monomial equations, when considered as a function of
$q$, is PORC. In fact I show that the number of solutions can be expressed in
the form $df(q)$ for some primitive polynomial $f(x)\in\mathbb{Z}[x]$, where%
\[
d=\alpha+\sum_{i=1}^{r}\alpha_{i}\gcd(q-n_{i},m_{i})
\]
for some rational numbers $\alpha,\alpha_{1},\alpha_{2},\ldots,\alpha_{r}$,
some integers $m_{1},m_{2},\ldots,m_{r}$ with $m_{i}>1$ for all $i$, and for
some integers $n_{i}$ with $0<n_{i}<m_{i}$ for all $i$. In addition I give an
algorithm for computing $d$ and $f$.

\section{Proof of Theorem 1}

Let $\varphi:\,$GL$(m,\mathbb{Q})\rightarrow\,$GL$(n,\mathbb{Q})$ give an
algebraic family of groups, and let $K$ be a finite field of order $q$, 
such that $\varphi _{K}(A)$ is defined for $A\in $\,GL$(m,K)$.
If we let $V$ be a vector space of dimension
$n$ over $K$ then we have a natural action of $\varphi _{K}(A)$ on $V$, and 
this gives an action of GL$(m,K)$ on $V$. We want to prove that the number
of orbits of GL$(m,K)$ on $V$, when considered as a function of $q$, is PORC.

The number of orbits is given by Burnside's Lemma. It is%
\[
\frac{1}{|\text{GL}(m,K)|}\sum_{A\in\,\text{GL}(m,K)}\text{fix}(\varphi
_{K}(A)),
\]
where fix$(\varphi_{K}(A))$ is $q^{d}$ where $d$ is the dimension of the
eigenspace of $\varphi_{K}(A)$ with eigenvalue 1. The number of orbits of
GL$(m,K)$ on $k$-dimensional subspaces of $V$ is given by the same formula,
where now fix$(\varphi_{K}(A))$ is the number of $k$-dimensional subspaces $W$
of $V$ such that $W\varphi_{K}(A)=W$. So we need to show that the functions
defined by these two formulae are PORC.

We simplify the problem as follows. There are only finitely many possible
types for matrices $A\in\,$GL$(m,K)$, and so it is sufficient to show that for
each type $T$%
\[
\sum_{A\text{ has type }T}\text{fix}(\varphi_{K}(A))
\]
is PORC. Actually, there is a slight problem here since that would only show
that the number of orbits had the form%
\[
\frac{f(q)}{|\text{GL}(m,q)|}%
\]
for some PORC function $f$. However, as Higman observes in \cite{higman60b},
if $k(x)$ is the quotient of two PORC functions, and if $k(x)$ only takes
integral values, then $k(x)$ is PORC. This is because a rational function of
$x$ which takes integral values for infinitely many integral values of $x$ is
a polynomial.

Since the size of the conjugacy class of an element of type $T$ in GL$(m,q)$
is a polynomial in $q$ which only depends on $T$, if we pick a set $S_{T}$ of
representatives for the conjugacy classes of type $T$ then it is only
necessary to show that%
\[
\sum_{A\in S_{T}}\text{fix}(\varphi_{K}(A))
\]
is PORC for each possible type $T$. As we saw in Section 6, if $A\in
\,$GL$(m,q)$ has type%
\[
T=\{(n_{1},S_{1}),(n_{2},S_{2}),\ldots,(n_{t},S_{t})\},
\]
then the conjugacy class of $A$ is determined by a set of eigenvalues
$\lambda_{1},\lambda_{2},\ldots,\lambda_{t}$ satisfying the equations%
\begin{equation}
\lambda_{i}^{q^{n_{i}}}=\lambda_{i}\;(i=1,2,\ldots,t),%
\end{equation}
and satisfying the non-equations%
\begin{equation}
\lambda_{i}^{q^{r}}\neq\lambda_{i}\;(0<r<n_{i},\;i=1,2,\ldots,t),%
\end{equation}
\begin{equation}
\lambda_{i}^{q^{r}}\neq\lambda_{j}^{q^{s}}\;(i\neq j,\;0\leq r<n_{i},\;0\leq s<n_{j}).%
\end{equation}
These eigenvalues can be taken to lie in $L=\,$GF$(q^{d})$ where $d$ is the least
common multiple of $\{n_{1},n_{2},\ldots,n_{t}\}$. As we run through all
possible solutions to these equations and non-equations in $L$ then the
conjugacy classes in GL$(m,K)$ determined by the solutions run through all
conjugacy classes of elements of type $T$ with each conjugacy class arising
the same number of times. (This number is given by equation (6) from Section
6.) So it is sufficient to show that%
\[
\sum_{\lambda_{1},\lambda_{2},\ldots,\lambda_{t}}\text{fix}(\varphi_{K}(A))
\]
is PORC where now the sum runs over all solutions in $L$ to the equations (7)
and non-equations (8), (9), and where $A\in\,$GL$(m,K)$ is chosen to lie in the
conjugacy class determined by the solution.

So let $\lambda_{1},\lambda_{2},\ldots,\lambda_{t}$ satisfy the equations (7)
and non-equations (8), (9). As we showed in Section 6, if $A$ lies in the conjugacy
class of GL$(m,K)$ determined by $\lambda_{1},\lambda_{2},\ldots,\lambda_{t}$
then $\varphi_{L}(A)$ has Jordan normal form%
\[
m_{1}J_{t_{1}}\oplus m_{2}J_{t_{2}}\oplus\ldots\oplus m_{s}J_{t_{s}}%
\]
for some positive integers $t_{1},t_{2},\ldots,t_{s}$, and some products
$m_{1},m_{2},\ldots,m_{s}$ of the form
\[
\lambda_{1}^{h_{1}}\lambda_{2}^{h_{2}}\ldots\lambda_{t}^{h_{t}}%
\]
where $h_{1},h_{2},\ldots,h_{t}$ are integer polynomials in $q$. The integers
$t_{1},t_{2},\ldots,t_{s}$ and the polynomials $h_{1},h_{2},\ldots,h_{t}$
depend only on the type $T$ and the characteristic of $K$.

Now consider the proof of Theorem 1 (a). In this case, for any given solution
$\lambda_{1},\lambda_{2},\ldots,\lambda_{t}$, fix$(\varphi_{K}(A))$ is $q^{d}$
where $d$ is the number of equations $m_{i}=1$ which are satisfied. The sequence
$m_{1},m_{2},\ldots,m_{s}$ and the size of the Jordan block associated with
each $m_i$ depend on the characteristic as well as on the type $T$,
so for the moment we assume that $K$ has fixed characteristic $p$. For every
subset $S \subseteq \{1,2,\ldots ,s\}$ let $U_S$ be the set of equations $m_i=1$
for $i\in S$ and let $V_S$ be the set of non-equations $m_i\neq 1$ for $i\notin S$.
Then Theorem 3 shows that the number of $\lambda_{1},\lambda_{2},\ldots,\lambda_{t}$ 
satisfying the equations and non-equations (7), (8), (9), $U_S$, $V_S$ is PORC 
when considered as a function of $q$. For all the solutions
fix$(\varphi_{K}(A))=q^{|S|}$. Every solution of (7), (8) and (9) satisfies
(7), (8), (9), $U_S$, $V_S$ for exactly one subset $S$, and so for each integer $d$ 
the number of solutions to (7), (8) and (9) for which fix$(\varphi_{K}(A))=q^{d}$ is 
PORC, and hence for each characteristic $p$ we obtain a PORC function $f_p(q)$ such that
\[
\sum_{\lambda_{1},\lambda_{2},\ldots,\lambda_{t}}\text{fix}(\varphi_{K}(A)) = f_p(q)
\]
whenever $q$ is a power of $p$. But as we saw in Section 5, there is a finite 
set of exceptional characteristics, and for all other characteristics the 
sequence $m_{1},m_{2},\ldots,m_{s}$ and the size of the Jordan block associated with
each $m_i$ depend only on $T$. So for each exceptional characteristic $p$ we 
obtain a PORC function $f_p(q)$ giving%
\[
\sum_{\lambda_{1},\lambda_{2},\ldots,\lambda_{t}}\text{fix}(\varphi_{K}(A))
\]
when $q$ is a power of $p$, and we obtain one further PORC function giving
this sum for all other characteristics. It follows that
\[
\sum_{\lambda_{1},\lambda_{2},\ldots,\lambda_{t}}\text{fix}(\varphi_{K}(A))
\]
is PORC as a function of $q$.

Finally consider the proof of Theorem 1 (b). Now fix$(\varphi_{K}(A))$ is the
number of $k$-dimensional subspaces $W$ of $V$ such that $W\varphi_{K}(A)=W$.
Eick and O'Brien \cite{eickobrien99} show that this number is given by a
polynomial in $q$, and that the polynomial only depends on the type of
$\varphi_{K}(A)$. Furthermore they give an algorithm for computing this
polynomial. As above, for the moment we assume that the characteristic
of $K$ is a fixed prime $p$. For any given solution $\lambda_{1}%
,\lambda_{2},\ldots,\lambda_{t}$, the type of $\varphi_{K}(A)$ is determined
by which equations $m_{i}=m_{j}$ hold (and which do not hold).

For every subset $S \subseteq \{(i,j)\,:\,1\leq i<j\leq s\}$ let $U_S$ be the set of 
equations $m_i=m_j$ for $(i,j)\in S$ and let $V_S$ be the set of non-equations 
$m_i\neq m_j$ for $(i,j)\notin S$.
Then Theorem 3 shows that the number of $\lambda_{1},\lambda_{2},\ldots,\lambda_{t}$ 
satisfying the equations and non-equations (7), (8), (9), $U_S$, $V_S$ is PORC 
when considered as a function of $q$. All the solutions to these equations
give matrices $\varphi_{K}(A)$ of the same type, and every solution to (7), (8) and 
(9) satisfies  (7), (8), (9), $U_S$, $V_S$ for exactly one subset $S$. So for every 
possible type $TT$ of $n\times n$ matrices, the number of solutions of (7), (8), 
(9) which give matrices $\varphi_{K}(A)$ of type $TT$ is PORC, and hence
for each characteristic $p$ we obtain a PORC function $f_p(q)$ such that
\[
\sum_{\lambda_{1},\lambda_{2},\ldots,\lambda_{t}}\text{fix}(\varphi_{K}(A)) = f_p(q)
\]
whenever $q$ is a power of $p$.

The rest of the proof of Theorem 1 (b) follows in the same way as the
proof of Theorem 1 (a).

\bigskip
It may help clarify the argument above if we look again at the example
given at the end of Section 5. We were looking at the Jordan form of
$A\otimes A$ when $A$ has type $T=\{(1,\{2\}),(1,\{3\})\}$. Matrices
$A=aJ_{2}\oplus bJ_{3}$ where $a^{q-1}=1$, $b^{q-1}=1$, $a\neq b$ give
a complete set of representatives for the conjugacy classes of matrices 
of type $T$ over GF$(q)$.

If $q=p^k$ where $p>3$ is prime then the Jordan form of the tensor square
of $aJ_{2}\oplus bJ_{3}$ is
\[
a^{2}J_{1}\oplus a^{2}J_{3}\oplus abJ_{2}\oplus abJ_{2}\oplus abJ_{4}\oplus
abJ_{4}\oplus b^{2}J_{1}\oplus b^{2}J_{3}\oplus b^{2}J_{5}.
\]
So in the notation used above $s=9$ and
\[
(m_1,m_2,\ldots ,m_9)=(a^2,a^2,ab,ab,ab,ab,b^2,b^2,b^2).
\]
Since the sequence $m_1,m_2,\ldots ,m_9$ has repetitions, many of the non-equations
$m_i\neq m_j$ are impossible. Also since $a\neq b$ many of the equations $m_i=m_j$
are impossible. But this makes no difference to the argument above since the
PORC formula giving the number of solutions to an impossible set of equations
and non-equations will be $0$. In our example $A\otimes A$ will have type
\[
\{(1,\{1,1,3,3,5\}),(1,\{2,2,4,4\})\}
\]
if $a^2=b^2$, and type
\[
\{(1,\{1,3\}),(1,\{2,2,4,4\}),(1,\{1,3,5\})\}
\]
if $a^2\neq b^2$.

Similarly, to determine how may of $m_1,m_2,\ldots ,m_9$ are equal to $1$,
we only need to determine which of $a^2,ab,b^2$ are equal to $1$. If $a^2=1$
or if $b^2=1$ then $ab=1$ is impossible, since $a\neq b$. So we need only
compute the number of solutions to $a^{q-1}=1$, $b^{q-1}=1$, $a\neq b$
when combined with each of the following five sets of equations and non-equations:
\[
a^2=1,\; ab\neq 1,\; b^2=1,
\]
\[
a^2=1,\; ab\neq 1,\; b^2\neq 1,
\]
\[
a^2\neq 1,\; ab\neq 1,\; b^2=1,
\]
\[
a^2\neq 1,\; ab=1,\; b^2\neq 1,
\]
\[
a^2\neq 1,\; ab\neq 1,\; b^2\neq 1.
\]

\bigskip
Christ Church

Oxford

OX1 1DP

michael.vaughan-lee@chch.ox.ac.uk


\begin{thebibliography}{99}                                                                                               
\bibitem {eickobrien99}B.~Eick and E.A. O'Brien, \emph{Enumerating $p$%
-groups}, J. Austral. Math. Soc.  Ser. A \textbf{67} (1999), 191--205.

\bibitem {eickwesche}B.~Eick and M.~Wesche, \emph{Enumeration of nilpotent
associative algebras of  class $2$ over arbitrary fields}, J. Algebra
\textbf{503} (2018), 573--589.

\bibitem {green55}J.A. Green, \emph{The characters of the finite general
linear groups}, Trans.  Amer. Math. Soc. \textbf{80} (1955), 402--447.

\bibitem {higman60}G.~Higman, \emph{Enumerating $p$-groups. {I}:
Inequalities}, Proc. London Math.  Soc. \textbf{(3) 10} (1960), 24--30.

\bibitem {higman60b}G.~Higman, \emph{Enumerating $p$-groups. {II}: Problems
whose solution is {PORC}},  Proc. London Math. Soc. \textbf{(3) 10} (1960), 566--582.

\bibitem {macdonald79}I.G. Macdonald, \emph{Symmetric functions and {H}all
polynomials}, The Clarendon  Press, Oxford University Press, New York, 1979.

\bibitem {newobvl}M.F. Newman, E.A. O'Brien, and M.R. Vaughan-Lee,
\emph{Groups and nilpotent  {Lie} rings whose order is the sixth power of a
prime}, J. Algebra \textbf{278}  (2004), 383--401.

\bibitem {vlee12a}Michael Vaughan-Lee, \emph{On {G}raham {H}igman's famous
{PORC} paper},  Internat. J. Group Theory \textbf{1} (2012), 65--79.

\bibitem {vlee13}Michael Vaughan-Lee, \emph{Enumerating algbras over a finite
field}, Internat. J. Group  Theory \textbf{2} (2013), 49--61.

\bibitem {vlee17}\mbox{Michael Vaughan-Lee, \emph{Choosing elements from finite
fields}}, arXiv.1707.09652 (2017).

\bibitem {witty}Brett Witty, \emph{Enumeration of groups of prime-power
order}, Phd thesis,  Australian National University, 2006.
\end{thebibliography}
\end{document}